\documentstyle[12pt]{article}
\makeatletter
\begin{document}
\newcommand{\genus}{\mathop{\mbox{genus}}}
\newcommand{\Endo}[1]{\mathop{\mbox{End}_{#1}}}
\newcommand{\rank}{\mathop {\mbox{rank}}\nolimits}
\newcommand{\tform}[1]{#1_{*}TM}
\newcommand{\tvttform}[1]{T^{*}M\otimes #1_{*}TM}
\newcommand{\one}[1]{ {\mbox{\boldmath $1$}}_{#1}}
\newcommand{\Gr}{\mathop {\mbox{\boldmath $Gr$}}\nolimits}
\newcommand{\re}{\mathop {\mbox{\boldmath $Re$}}\nolimits}
\newcommand{\im}{\mathop {\mbox{\boldmath $Im$}}\nolimits}
\newcommand{\QED}{$\Box$ \\}
\newcommand{\tildef}{\tilde{f}}
\newcommand{\Ztwo}{\mbox{\boldmath $Z_{2}$}}
\newcommand{\R}{\mbox{\boldmath $R$}}
\newcommand{\C}{\mbox{\boldmath $C$}}
\newcommand{\bfZ}{\mbox{\boldmath Z}}
\newcommand{\qH}{\mathop {\mbox{\boldmath $H$}}\nolimits}
\newcommand{\bfQ}{\mbox{\boldmath Q}}
\newcommand{\bfN}{\mbox{\boldmath N}}
\newcommand{\bfP}{\mbox{\boldmath P}}
\newcommand{\bfT}{\mbox{\boldmath T}}
\newcommand{\sqdz}{\sqrt{dz}}
\newcommand{\sqdzb}{\sqrt{d\bar{z}}}
\newcommand{\hD}{{\cal{D}}^\frac{1}{2}}
\newcommand{\aiso}{\mathop {\cal{A}} \nolimits}
\newcommand{\calD}{{\cal{D}}}
\newcommand{\imH}{\mathop{\rm{im}}(\qH)}
\newcommand{\idx}{\mathop {\rm ind}\nolimits}
\newcommand{\setB}[1]{\left\{ #1 \right\}}
\newcommand{\divg}{\mathop {\rm div}\nolimits}
\newcommand{\cS}{{\cal S}} 
\newcommand{\fref}[1]{Figure~\ref{#1}}
\newcommand{\cM}{{\cal M}}
\newcommand{\MF}{{\cal MF}}
\newcommand{\mref}[1]{{\rm (\ref{#1})}}
\newcommand{\cT}{{\cal T}}
\newcommand{\s}{\mathop {\rm R}\nolimits}
\newcommand{\Sg}{(\Sigma, g)}
\newcommand{\RS}{\R_{+}^{\cal S}}
\newcommand{\Dif}{\mathop{\rm Diff_0}}
\newcommand{\pone}{\psi_{1}}
\newcommand{\ptwo}{\psi_{2}}
\newcommand{\None}{N_{1}}
\newcommand{\Ntwo}{N_{2}}
\newcommand{\Sione}{\Sigma_{1}}
\newcommand{\Sitwo}{\Sigma_{2}}
\newcommand{\II}{\mathop {\rm I\!I}\nolimits}
\newcommand{\IIone}{\mathop {\rm I\!I}_{1}\nolimits}
\newcommand{\IItwo}{\mathop {\rm I\!I}_{2}\nolimits}
\newcommand{\Sig}{(\Sigma, g)}
\newcommand{\rh}{\sim_{r}}

\newcommand{\n}{\noindent}
\newtheorem{obs}{Observation}
\newtheorem{theorem}{Theorem}
\newtheorem{assertion}{Assertion}
\newtheorem{proposition}{Proposition}
\newtheorem{Rem}{Remark}
\newenvironment{remark}{\Rem \rm}{\endRem}
\newtheorem{Conv}{Convention}
\newenvironment{convention}{\Conv \rm}{\endConv}

\newtheorem{lemma}{Lemma}
\newtheorem{definition}{Definition}
\newtheorem{claim}{Claim}[section]
\newtheorem{corollary}{Corollary}[section]
\newtheorem{coro}{Corollary}[section]
\newtheorem{observation}{Observation}[section]
\newtheorem{conjecture}{Conjecture}[section]
\newtheorem{question}{Question}[section]
\newtheorem{example}{Example}[section]
\begin{center}
{\large{\bf QUADRATIC DIFFERENTIALS, QUATERNIONIC FORMS, AND SURFACES} }
\vskip.25in
George I. Kamberov \\ [2mm]
\footnotesize{
Department of Mathematics \\
Washington University \\
Saint Louis, MO 63130} 
\end{center}

\vskip .25in

\section{Introduction} \label{intro} 
In this paper we address three problems stemming from the fundamental question: 
{\em Let $M$ be an oriented surface, how much data do we need to identify an immersion 
$f: M \rightarrow \R^3$ or its shape?} 

A generic immersion $f$ is determined up to homothety and translation by its
conformal class and its tangent plane map $T_f : M\rightarrow
\Gr(2,3)$ whose value $T_{f}(p)$ at each point $p\in M$ is the plane
tangent to the image surface $f(M)$. There are exceptions, called {\bf
  Christoffel immersions}. Their classification is known as 
{\bf Christoffel's problem}. The problem was posed and to a large extent 
solved locally in \cite{christoffel}. A short complete solution of Christoffel's
problem, including closed surfaces, is presented in
Section~\ref{christ}.

Bonnet noticed in \cite{Bonnet} that if two immersions induce the same 
first fundamental form and the same mean curvature function then they
are congruent unless they are quite exceptional. {\bf Bonnet's
  problem} is to classify all such exceptional immersions.  A local
classification of the umbilic-free and the constant mean curvature
exceptions follows from \cite{Bonnet, Cartan-pr-curv,
  Chern-Deform-Pr-Curv}. Until now it was not known whether there
exist other exceptions. The existence and description of the new examples is 
discussed in Section~\ref{sect:bonnet}.

Christoffel's and Bonnet's problems inevitably lead to the study of
{\bf isothermic immersions} \cite{darboux, pipeka}.  An umbilic-free
immersion is called umbilic-free classical isothermic immersion if
every point $p\in M$ admits isothermal local coordinates in which its
second fundamental form is diagonal. There exist classical definitions
of isothermic immersions with umbilic points but no definitive one. Up
to now it was not clear how to motivate and define global isothermic
immersions.  In this paper we use a hitherto unsuspected connection
between isothermic immersions and quadratic differentials to study
global isothermic immersions and immersions with umbilic points,
branch points, and ends (see Section~\ref{sec:is}). The presented
theory is designed to solve the targeted geometric problems and to
resolve some quirks in the folklore surrounding the classical
isothermic immersions.

The main tool used in this paper is the quaternionic calculus on
surfaces \cite{pipeka}.  The necessary material is summarized in
Section~\ref{sec:afqd}.  Some of the results in this paper were
announced in \cite{hqd}. 

Unless explicitly specified otherwise all immersions in this 
paper are assumed to be at least $C^{3}$. 
On the other hand many of the constructions and arguments 
do not need much smoothness.  Weaker smoothness requirements, regularity, 
and partial regularity are discussed in \cite{sqd}.

\section{Quaternionic Calculus, Tensors, and Function Theory}
\label{sec:afqd}
A brief summary of the quaternionic calculus introduced 
in \cite{pipeka} is given below. (Very detailed treatments were 
given in \cite{ulrichsnB, saopaulo}). 

Throughout this paper $M$ will denote a Riemann surface with a complex structure $J$. 
Let $f$ be a conformal immersion of $M$ into $\R^3 = \imH$, and let $N$ be the
oriented unit normal vector field of $f$,  a quaternionic valued
one-form $\kappa$ on $M$ is called conformal (anti-conformal) with
respect to $f$, if $*\kappa = N\kappa$ ($*\kappa = -N\kappa$), where 
$*\kappa := \kappa\circ J$. More
generally, an   $\qH$ valued one-form $\varphi$ defined on a
subset of a Riemann surface is called conformal if there exists 
a quaternionic valued function $N$ such that $*\varphi = N\varphi$. 
Given a conformal immersion, then every quaternionic valued
one-form $\tau$ is the sum of a conformal  
part and a anti-conformal part. In particular, let $N$ be the Gauss
map of the immersion $f$ then the form $dN$ is tangential valued and 
\begin{equation} \label{eqn:Weingarten}
dN = -Hdf + \omega, 
\end{equation}
where $H$ is the usual mean curvature function of $f$ and $\omega$ is 
the anti-conformal part of $dN$. 
The anti-conformal part $\omega$ of $dN$ is called the {\em Hopf form} of
the immersion. See Remark \ref{hopf} below. 

Let $f$ be a conformal immersion with normal $N$. Every quaternionic 
valued one-form $\tau$ has a splitting $\tau =\tau^{\top} + \tau^{\perp}$ 
into a {\em tangential} and {\em transversal} part with respect to $f$. 
Note that a form is tangential if an only if $N\tau = - \tau N$,  while 
$\tau^{\top} = 0$ if and only if $N\tau = \tau N$. We denote the bundle 
of tangential one-forms by $T^{*}M\otimes \tform{f}$. 

Let $\tau$ be a tangential valued differential 
one-form. We denote the 
transversal and the tangential components  of $d\tau$, with respect to $f$,
by $(d\tau)^{\perp}$ and $(d\tau)^{\top}$, respectively.
\begin{observation} \label{dtau} Let $f$ be a conformal immersion, then for 
every tangential valued differential form $\tau$ we have:
\begin{equation} \label{dtauperp} 
(d\tau)^\perp = \left\{\begin{array}{ll} 
       \frac{1}{2}(\omega\wedge \tau -\tau\wedge\omega)N & 
                                 \mbox{if $\tau$ is anti-conformal} \\
        \\

        \frac{H}{2}(\tau\wedge df - df\wedge \tau)N & 
                                 \mbox{if $\tau$ is conformal}
                          \end{array}
                     \right.
\end{equation}
where $H$ and $\omega$ are the mean curvature and the Hopf form of the immersion $f$, 
respectively. 
\end{observation}
{\bf Proof} \,  \begin{eqnarray} \label{eqn:perp1}
(d\tau)^\perp &=& -\frac{1}{2}\left((d\tau)N +
                   \overline{(d\tau)N}\right)N \nonumber \\
 &=& -\frac{1}{2}\left((d\tau)N + Nd\tau \right)N. 
\end{eqnarray} 
The form $\tau$ is tangential. Differentiating $N\tau + \tau N = 0$  
yields $(d\tau)N + Nd\tau  =  \tau \wedge (dN) -   (dN) \wedge\tau$. 
Using this identity  and \mref{eqn:Weingarten} we obtain 
\[
-\frac{1}{2}\left((d\tau)N + Nd\tau\right) = \left\{\begin{array}{ll} 
\frac{1}{2}(\omega\wedge \tau -\tau\wedge\omega) & 
                                 \mbox{if $\tau$ is anti-conformal} \\
        \\
\frac{H}{2}(\tau\wedge df - df\wedge \tau) & 
                                 \mbox{if $\tau$ is conformal}
                          \end{array}
                     \right.
\]
Substituting back in \mref{eqn:perp1} completes the proof. 
\QED
\begin{definition}\label{def:spin}
Two conformal immersions $f,\tilde{f}:M\to\R^3 = \imH$ are called {\em spin-equivalent}
if there exists $\lambda:M\to\qH^{*}$ so that $d\tilde{f}=\bar{\lambda}d f\lambda$.
We call $\tilde{f}$ a {\em spin-transform} of $f$.
\end{definition}
Note that if $f$ and $\tilde{f}$  are  in the same regular homotopy
class then they are spin equivalent. In particular, $\pi_{1}(M) = 0$ then 
every two conformal immersions are spin equivalent.  

\begin{remark} \label{rem:int} The immersion $f$ defines an isomorphism 
between $\tvttform{f}$ and $T^{*}M\otimes T^{*}M$, sending 
a tangential valued form $\tau$ to the covariant tensor $<
\tau(\cdot)|df(\cdot)>$. (Here $<\cdot \, | \, \cdot >$ is the Euclidean 
scalar product in $\R^3$.)

A form $\kappa\in T^{*}M\otimes \tform{f}$  is conformal if and only
if it is a linear combination of $df$ and $*df$;  
the tensor $< \tau(\cdot)|df(\cdot)>$ is antisymmetric if and only 
if $\tau = a*df$ for some real $a$. 
The tensor $< \tau(\cdot)|df(\cdot)>$ is symmetric trace-free  
if and only if $\tau$ is anti-conformal with respect to $f$. Moreover, 
$< \tau(\cdot)|df(\cdot)>$ is divergence free (with respect to 
the metric induced by $f$) if and only if $(d\tau)^{\top}=0$. 

Suppose that $\tilde{f}$ is obtained from $f$ by a spin transform, thus
  $d\tilde{f} = \overline{\lambda}df\lambda$. Let $\tilde{\tau}$ be a 
  and $\tau$ be two $\qH$-valued one forms on $M$, tangential to 
$\tilde{f}$ and $f$ respectively, then 
\[
<\tilde{\tau}(\cdot)|d\tilde{f}(\cdot)> = <\tau(\cdot)|df(\cdot)>
\mbox{if and only if $\tau = \lambda \tilde{\tau} \overline{\lambda}$}.
\]
\end{remark} 
Remark~\ref{rem:int} indicates that the tangential
anti-conformal forms are an extrinsic geometric representation of the
quadratic differentials, and $(d\ )^{\top}$ is the extrinsic operator
$\overline{\partial}$. This analogy can be taken further.  In the
context of the intrinsic conformal geometry on a Riemann surface $M$
all complex valued objects take values in a fixed complex plane, $\C$.
On the other hand a conformal immersion $f$ of $M$ into $\R^3$ defines
two fields of complex planes, the field $f_{*}TM$ of tangent planes,
and the field of normal planes, $\R\oplus\R \, N$. The intrinsic
conformal objects can be represented as $f_{*}TM$-valued sections or
as $\R \oplus \R N$-valued sections. For example, let $w=a + ib$ be a
complex valued function on the Riemann surface $M$. Given a conformal
immersion $f$ with normal $N$, we can consider this function $w$ as
the section $(a + bN)df$ of the bundle of tangential and conformal
forms. 
\begin{observation} \label{obs:holof}
If $f$ is a conformal immersion of $M$ into 
$\R^3$  then a complex valued function $w = a + ib$ is holomorphic if and only if 
$\left(d(a + bN)df\right)^{\top} = 0.$
\end{observation}
{\bf Proof} \,  This is a local statement. Let $z = x + iy$ is a holomorphic coordinate chart on $M$,
then using \mref{eqn:Weingarten} we obtain 
\[
(d(a + bN)df)^{\top}(\partial_{x}, \partial_{y}) = 
-\left(\frac{\partial a}{\partial y} + \frac{\partial  b}{\partial x}\right)
\frac{\partial f}{\partial x} + 
\left( \frac{\partial a}{\partial x} - \frac{\partial  b}{\partial y}\right)
\frac{\partial f}{\partial y}.
\]
Therefore, $(d(a + bN)df)^{\top} = 0$ if and only if $a + ib$ is
holomorphic. \QED 

The following observation is based on ideas introduced in \cite{ulrichsnB}. 
 (See \cite{ulrichsnB, saopaulo} for proofs.) 
The bundle of {\em normal-valued quadratic
  differentials} 
associated with $f$ is the bundle of fiber-preserving maps defined by:
\[\R\oplus\R N^{2,0} := \left\{\kappa :TM\rightarrow \R\oplus\R \! N \, | \, 
\kappa(aX + bJX) = \kappa(X)(a + bN)^{2}, \, \forall X\in TM \right\}\]
The immersion $f$ defines a bundle isomorphism 
$\nu_f : \R\oplus\R \, N \rightarrow M\times \C$, by $\nu_f(a+bN) = a+ib$. 
The map sending $\kappa\in \R\oplus\R N^{2,0}$ to $\nu_{f}\circ \kappa$ is 
an isomorphism between $\R\oplus\R N^{2,0}$ and $T^{2,0}M^{*}$. 
A section $\kappa\in \Gamma(\R\oplus\R N^{2,0})$ is called a holomorphic 
normal-valued quadratic differential if the 
quadratic differential $\nu_{f}\circ\kappa$ 
is holomorphic. Once a conformal immersion is chosen it is convenient to
abuse notation and call normal-valued quadratic differentials simply 
quadratic differentials. Moreover, we have
\begin{observation} Every  conformal immersion, $f$, induces an 
isomorphism between the bundle of anti-conformal tangential one forms and 
the bundle of (normal) quadratic differentials, 
sending the one-form $\tau$ to the quadratic differential $df \tau$. 

Furthermore a quadratic
  differential $\kappa$ is holomorphic if and only if the tangential 
part of $d(df\backslash\kappa)$ vanishes, that is if and only if 
$\mbox{$(d(df\backslash\kappa))^{\top}=0$}$. 
\end{observation} 
\begin{remark} \label{hopf} Given a complex valued quadratic differential 
$q$ and a conformal immersion $f$ we will use the notation $df\backslash q$ 
for the tangential valued one form $df\backslash\nu_{f}^{-1}\circ q$. 
Furthermore notice that for every anti-conformal tangential form $\tau$ the 
two-zero part of the symmetric trace-free tensor $<\tau(\cdot)|df(\cdot)>$ 
is precisely the quadratic differential $\nu_{f}\circ df \tau$. 
The  isomorphism $q \rightarrow df\backslash q$ identifies the classical Hopf
  differential of $f$ with (minus) the Hopf form of $f$. Indeed, the Hopf differential 
of a conformal immersion $f$ is the two-zero part, $\II^{2,0}$,  of the second fundamental 
form $\II = -dN$. On the other hand \mref{eqn:Weingarten} implies that the tensor 
$<-\omega(\cdot)|df(\cdot)>$ is the trace-free part of $\II$ and so the Hopf 
differential of $f$ is $\II^{2,0} = -\nu_{f}\circ df \, \omega$, that is 
$\omega = -df \backslash \II^{2,0}$. 
\end{remark}
\section{Isothermic Surfaces}
\label{sec:is}
\begin{definition} \label{isoextended} A conformal immersion $f: M \rightarrow \imH$ is 
called an {\bf extended globally 
  isothermic immersion} if there exists a subset $S\subset M$ and a
map $f^* : M\setminus S \rightarrow \imH$  such that $M\setminus S$ is
dense in $M$,  
$df^* \not\equiv 0$,  and $df df^*$ is a holomorphic quadratic differential 
on $M\setminus S$. The map $f^*$ is called {\bf dual} to the immersion $f$, 
and the set $S$ is called the {\bf singular set} of $f^{*}$. 
\end{definition} 
\begin{remark} \label{rem:branch} 
Note that the dual $f^*$ of $f$, satisfies 
$df^{*} = df\backslash df df^{*}$ on $M' := M\setminus S$, and so $df^{*}$ is anti-conformal 
tangential with respect to the restriction of $f$ on $M'$. Furthermore, 
the zero set $Z = \setB{p\in M' | df^{*}_{p} = 0}$ is precisely the 
zero locus of the holomorphic quadratic differential $df df^*$ and so $Z$ is 
discrete, if not empty. Thus the dual map $f^{*}$ defines a branched 
conformal immersion of the Riemann surface $M'$. The branch locus of 
$f^{*}$ is precisely the zero set $Z$,  and the oriented normal of the immersion 
$f^{*}$ is $N^{*} = -N$. 
\end{remark}
\begin{remark} \label{singularset} An immersion $f$ is 
  extended globally isothermic if and only if there exists a holomorphic
  quadratic differential $q\not\equiv 0$, defined on a dense subset $M' =
  M\setminus S$, such that $df\backslash q$ is exact.
  Definition~\ref{isoextended} can be generalized. One could
  require that $df\backslash q$ is closed instead of exact but presently this
  appears justified only in the study of pathologies, for example immersions
  with uncountably many ends.  Note also that both the immersion $f$ and its
  dual $f^{*}$ may have only first order derivatives, and in particular, they
  may not have well defined second fundamental forms, curvature lines, and
  umbilic points in the usual sense.  Generalized isothermic immersions and
  the related regularity theory are studied in \cite{sqd} in connection with
  weak isometric immersions and the realization of conformal structures.
  Smooth local (umbilic-free) isothermic immersions in co-dimension two are studied in
  \cite{franzudo}.
\end{remark}

Let $\overline{\R^{3}} = \R^{3} \cup \setB{\infty}$ be the one-point compactification of 
$\R^{3}$ with the natural topology. The following observation is now trivial: 
\begin{observation} Let $f: M \rightarrow \imH$ be a conformal 
immersion such that there exists a set $S\subset M$ and a 
map $f^{*}: M\setminus S \rightarrow \R^{3}$, such that 
$df df^{*}$ is a nontrivial holomorphic quadratic
differential on $M\setminus S$, which extends to a meromorphic 
quadratic differential on $M$ whose set of poles is $S$. Then $f$ is an 
extended isothermic immersion with dual $f^{*}$.

 The map $f^{*}$ is a branched immersion of $M$ with ends. 
Namely, $f^{*} : M \rightarrow \overline{\R^{3}}$ is continuous and 
$(f^{*})^{-1}(\infty) = S$ and so the poles of $df df^{*}$ are precisely the ends 
of $f^{*}$. The set $Z$ of zeroes of the differential 
$df df^{*}$ is the set of branch points of the extended immersion
$f^{*}$. Finally, the branch points and the ends of $f^{*}$ are
umbilic points of $f$ but $f$ may have other umbilic points.  
\end{observation} 
Every rotationally symmetric immersion $f_{s}$ is  extended globally isothermic:
indeed, $f_{s}$ has a rotationally symmetric dual $f^{*}_s$. In particular, 
every ellipsoid  of revolution $f_{s}$ is extended globally 
isothermic: the dual $f^{*}_{s}$ has 
precisely two ends according to the nomenclature adopted 
in this paper. General genus zero immersions including non-rotational 
ellipsoids are discussed in \cite{sqd}. 
\begin{definition} \label{iso} A  
conformal immersion $f$ is called {\bf globally   isothermic}  if there exists a 
nontrivial holomorphic quadratic differential $q$ on $M$ such that $df\backslash q$ is exact. 
\end{definition}
The class of global isothermic immersions is precisely the class used in \cite{pipeka}, 
see Remark~\ref{rem:branch}. In the rest of this paper 
{\em isothermic immersion} is used as a synonym of {\em global isothermic
  immersion}. 

The defining property of classical  isothermic 
immersions is that in a neighborhood of every non-umbilic point 
they admit an isothermal coordinate system in which the second fundamental form 
is diagonal, equivalently  in a neighborhood of a every non-umbilic point 
one can represent the Hopf differential of the immersion as a real 
multiple of a holomorphic quadratic differential. 

Let $f$ be a global isothermic immersion with dual  $f^*$, and a Hopf form $\omega$. 
Thus from $d(df^{*})=0$ and from \mref{dtauperp} it follows 
that 
\begin{equation} \label{d(df^{*})}
df^*\wedge \omega - \omega\wedge df^* = 0.
\end{equation} 
Away from the branch locus $Z$ of $f^*$, the bundle of tangential anti-conformal forms with respect 
to $f$ is spanned by $df^{*}$ and $*df^*$. Thus  $\omega = adf^* + b*df^{*}$ where $a$ and $b$ 
are real-valued and smooth away form the branch locus $Z$. From \mref{d(df^{*})} 
it follows that $b = 0$ and so $\omega = adf^*$ and so the Hopf differential
of $f$ equals minus $a \nu_{f}\circ  df df^{*}$. The classical formula $dN = H^{*} df^{*} - H df$, 
where $H^{*}$ is the mean curvature of $f^{*}$, remains valid 
away from the branch locus of $f^*$. Therefore, $a = H^{*}$, away 
from the branch points of $f^{*}$. 

In a neighborhood of $p\in M\setminus Z$ choose isothermal 
coordinates $z=x+iy$, with respect to the metric induced by $f$, in which the quadratic 
differential $q = \nu_{f} \circ df df^{*} = dz^{2}$. In this coordinate 
system the Hopf differential $\II^{2,0}$ of $f$ equals $adz^{2}$. 
The coefficient $a$ is real valued. Therefore, the second fundamental form of 
$f$ is diagonal with respect to the $(x,y)$ coordinates. Note that $p$ could be 
an umbilic point for $f$. Thus the following observation is true for global isothermic immersions: 
\begin{observation} \label{obs:hqd=>}
  Every global isothermic immersion $f$ defines a global holomorphic
  quadratic differential $df df^{*}$, where $f^{*}$ is dual to $f$.  The Hopf
  differential of $f$ is a real multiple of $df df^{*}$ (smooth multiple on $M\setminus Z$). For every point $p\in M\setminus Z$ there exists an
  isothermal coordinate system in which the second fundamental form of $f$ is
  diagonal.
\end{observation}
This observation indicates that while in the classical theory umbilics are avoided and 
excluded, they are generically quite benign. The only {\em interesting} umbilic points are the 
branch points of the dual. 

Observation~\ref{obs:hqd=>} has a local inverse. The data determining an isothermic immersion 
consists of the tangent map of the immersion along a curve, which is 
not a curvature line, together with a holomorphic quadratic differential 
encoding the properties of the lines of curvature of the immersion. Recall that 
with every quadratic differential $q$ on a Riemann surface $M$ we associate 
two mutually orthogonal line foliations (possibly singular). These foliations are 
basic objects in Teichm\"{u}ller theory, they are called the {\em principal stretch}  
foliations of $q$. We will call a smooth curve ${\cal C}$ {\bf non-characteristic} for 
$q$ if it is transversal to the principal stretch foliations of $q$. 
We prove the basic existence result in the analytic category below, and indicate the 
general result in Remark~\ref{rem:gen:exist} below. 
\begin{theorem}\label{obs:hqd<=}
Given a holomorphic quadratic differential $q$ on a 
Riemann surface $M$, let ${\cal C}$ be a real analytic curve in 
$M\setminus Z$ which is non-characteristic for $q$, then for every  
real analytic, conformal,  rank two 
form $\varphi_0$ defined along ${\cal C}$, there exists a unique, up 
to translation in space, isothermic immersion $\tilde{f}$ 
defined in a neighborhood of ${\cal C}$ such that $d\tilde{f}\backslash q$ is 
the differential of the dual of $\tilde{f}$ and 
$d\tilde{f} = \varphi_0$ along ${\cal C}$.  
\end{theorem}
{\bf Proof} \,  It suffices to find an
immersion $\tilde{f}$ such that $d\tilde{f}\backslash q$ is exact; this is
equivalent to showing that the normal component of $d(d\tilde{f}\backslash q)$
vanishes.  To find $\tilde{f}$ we will solve
an initial value problem with initial data $d\tilde{f} = \varphi_0$
along ${\cal C}$.  Choose a background real analytic conformal
immersion $f$, such that $df = \varphi_0$ along ${\cal{C}}$, and let 
$\tau = df\backslash q$. 
  Then we search
 for $\tilde{f}$ as a spin transform of
$f$. Equivalently we need to find a quaternionic valued function
$\lambda$ such that $d\tilde{f} = \overline{\lambda} df \lambda$,  
$d\tildef \backslash q = \lambda^{-1} (df\backslash q) 
\overline{\lambda^{-1}}$ is exact, 
and $\lambda = 1$ along ${\cal C}$. 
The function $\lambda$ must satisfy the differential equations 
\begin{eqnarray*}
  d (\overline{\lambda} df \lambda) &=& 0 \\
  d ( \lambda^{-1} \tau \overline{\lambda^{-1}})^{\perp} &=& 0.
\end{eqnarray*}
Thus the system for $\lambda$ is 
\begin{eqnarray}  \label{eq:1.a}
  \im\left(\overline{\lambda} df \wedge d\lambda\right) &=& 0 \\
  \re\left((d\lambda^{-1}\tau\overline{\lambda^{-1}})\lambda^{-1}N
\lambda\right) &=& 0, 
\label{eq:1.b}
\end{eqnarray}
A direct computation shows that equation \mref{eq:1.b} is equivalent to
the equation
\begin{equation}
  \label{eq:1.b.2}
  \re\left(\im\left((d\lambda)\lambda^{-1}\wedge \tau\right)N\right) + 
(\omega\wedge\tau - \tau\wedge\omega)/4 = 0
\end{equation}
where $\omega$ is the Hopf form of $f$. Thus we must solve  
\mref{eq:1.a}, \mref{eq:1.b.2}, with Cauchy data $\lambda = 1$, 
along ${\cal C}$. In local coordinates $z = x + iy$ the differential system is 
\begin{eqnarray}
  \label{eqna:coords:1}
  \im\left(\overline{\lambda}\left(
 \frac{\partial f}{\partial x} \frac{\partial \lambda}{\partial y} - 
 \frac{\partial f}{\partial y} \frac{\partial \lambda}{\partial x}
\right)
\right) &=& 0 \\
\re\left(\im\left(
\frac{\partial \lambda}{\partial x}\lambda^{-1} \tau\left(\frac{\partial }{\partial y}\right) - 
 \frac{\partial \lambda}{\partial y}\lambda^{-1} \tau\left(\frac{\partial }{\partial x}\right)
\right)N
\right) + \cdots &=& 0 \label{eqna:coords:2}
\end{eqnarray}
We can use \mref{eqna:coords:1} and \mref{eqna:coords:2} to compute the principal 
symbol of the system of differential equations
\mref{eq:1.a}, \mref{eq:1.b.2} at $\lambda = 1$, and a  
covector $\xi$.  Indeed, let $\xi = \xi_{1} dx + \xi_{2} dy$ then the symbol is the map 
\[
\sigma(1,\xi) : \qH \rightarrow \qH
\]
defined by 
\begin{equation} \label{eqn:symb:coor} 
\sigma(1,\xi)(\alpha) = \re\left(\im\left(\alpha B(\xi)\right)N\right) + 
\im\left(A(\xi)\alpha\right), 
\end{equation}
where 
\[
A(\xi) = \xi_{2}\frac{\partial f}{\partial x} - \xi_{1}\frac{\partial f}{\partial y}
\hspace{1em} \mbox{and} \hspace{1em}
B(\xi) = \xi_{1}\tau\left(\frac{\partial }{\partial y}\right) 
      -  \xi_{2}\tau\left(\frac{\partial }{\partial x}\right). 
\]
Note that both $A(\xi)$ and $B(\xi)$ are tangential to 
the immersion and in particular they are $\imH$ valued, and 
so $\sigma(1,\xi)(\alpha)=0$ precisely when $\alpha$ is imaginary 
valued and is collinear to both $B(\xi)$ and $A(\xi)$. This implies 
that the vector dual to $\xi$, with respect to the first 
fundamental form induced by $f$, is tangential to the  principal stretch 
foliation of the quadratic differential $q$. Thus the characteristic variety 
of the differential system for $\lambda$ consists of two orthogonal lines, 
each one tangent to the stretch directions of $q$. Therefore, the initial value problem 
with Cauchy curve transversal to the principal stretch foliations of $q$ is well posed. 
\QED

In summary the space of all isothermic immersions generated by $q$
is parameterized by four functions of one real variable. These
functions represent the differential of the immersion along the
initial curve. (See \cite{saopaulo} and \cite{sqd} for more details.)
Compare with the result by Weingarten see
\cite{darboux}, page 262, and the references there.

\begin{remark} \label{rem:gen:exist}
In fact one can remove the technical assumption that the 
curve ${\cal C}$ omits the zero locus $Z$. After modifying the definition 
of non-characteristic curve one can still obtain an existence result. 
Furthermore, the existence results are true in the smooth category. 
The proofs in these  cases require a whole new set of techniques. 
(See \cite{sqd}.)
\end{remark} 
\section{Christoffel's Problem} \label{christ}
\begin{theorem} \label{christsth} Every Christoffel immersion $f$ is 
globally isothermic. Vice  versa let $f$ be a globally isothermic immersion 
with dual $f^{*}$ and let
  $Z$ be the branch locus of $f^{*}$, then either $Z$ is empty or it
  is a set of isolated points and every point in $Z$ is umbilic for
  $f$. Furthermore, $f$ is a  Christoffel immersion on 
  $M\setminus Z$. 
\end{theorem} 
{\bf Proof} \,  Suppose $f$ and $\tilde{f}$ are 
conformal immersions of the same Riemann surface into $\R^3$ 
such that $T_f = T_{\tildef}$. In particular, $d\tildef$ is tangential,  
and conformal or anti-conformal with respect to $f$. If $d\tilde{f}$ is 
anti-conformal with respect to $f$ then $f$ is isothermic and
$\tilde{f}$ is a dual to $f$. To complete the proof we need to
consider the case when $d\tilde{f}$ is conformal 
with respect of $f$. In this case the normal to $\tilde{f}$ is $N$ and
there exist real 
valued functions $a$ and $b$ such that 
\[d\tilde{f} = (a + bN)df.\] 
Clearly $(d d\tildef)=0$. The tangential part of this identity implies 
that $u=a + ib$ is a holomorphic function defined on $M$. (See 
Observation~\ref{obs:holof}.) On the other hand  from 
\mref{dtauperp}, it follows that
\[0 = (d d\tildef)^\perp = H 
                          \left(d\tildef\wedge df - df\wedge d\tildef 
                          \right)/2, \] 
where $H$ is the mean curvature of $f$. 
If $H_p \neq 0$ at some point $p\in M$, then 
\[0=df\wedge d\tildef - d\tildef\wedge df = b(df\wedge *df - *df\wedge df)\]  
in a neighborhood of this point, and so $u$ is constant and real in this 
neighborhood and hence $u$ must be constant and real on the whole surface $M$.
This implies that $\tildef$ is just a real scaling of $f$. This is impossible 
if $f$ and $\tildef$ are a Christoffel pair. Thus $f$ must be minimal, 
but every minimal immersion is isothermic.

The proof of the second part of the theorem follows immediately from 
Remark~\ref{rem:branch}. 
\QED

Recently Holy Bernstein and Gary Jensen  used complex frames (\cite{gary}) 
to  give another proof of the first part of this theorem \cite{holy}. 

Let $M$ be a closed surface. If $f$ is a Christoffel immersion suppose that $\tilde{f}$ 
is another immersion such that $T_{f} = T_{\tilde{f}}$ but $f$ is not
obtained from $\tilde{f}$ by homothety and translation, then $f$ must
be isothermic, $\tilde{f}$ must be dual to $f$, and so $df d\tilde{f}$
must be holomorphic. Every holomorphic quadratic differential on a
closed surface whose genus is not equal to one must vanish somewhere.
Therefore, if $M$ is closed and $\genus(M)\neq 1$, then $\tilde{f}$ must
be a branched immersion, and not an immersion. Thus we obtain:
\begin{corollary} \label{coro:compact:christ} Let $M$ be a compact oriented surface whose genus is 
  not one, then every immersion $f$ of $M$ into $\R^3$ is
  determined uniquely up to homothety and translation by its conformal
  class and its tangent planes map, $T_f$.
\end{corollary}
\begin{remark} \label{tori} This result is false for tori. The 
obvious counterexamples are CMC tori. 
\end{remark}
\section{Bonnet's problem} \label{sect:bonnet} Two non-congruent 
conformal immersions $f_{\pm}$ are called Bonnet mates if they induce
the same metric and have the same mean curvature $H_{+}=H_{-}$. 
A complete local classification of the umbilic-free 
immersions which admit Bonnet mates was obtained in \cite{Bonnet,
Cartan-pr-curv, Chern-Deform-Pr-Curv}. The construction of 
explicit examples is often a difficult task (see \cite{bobe, rou}).
Some global results were obtained in 
\cite{trib_and_lawson, ARoss, kamberov-bonnet-96, eradgga}. 
It is conjectured that there are no compact Bonnet mates. 
A first step in the verification of this conjecture is 
to study Bonnet mates with umbilic points. 
The methods used by Bonnet, Cartan, and Chern can not work for
surfaces with umbilics. Moreover, until now it was not known whether
there exist non-constant mean curvature Bonnet mates with umbilic
points. In particular it was not known whether the Bonnet-Cartan-Chern
classification includes all the possible Bonnet mates. 
A general construction which yields all immersions admitting 
Bonnet mates -- with or with out umbilics was presented in \cite{pipeka}. 
A natural question is: Does this construction yield new hitherto unknown
examples? 

Let $f_{\pm}$ be Bonnet mates defined on the surface $M$, then 
up to a possible rotation of one of the immersions, say $f_{+}$,  we have 
\begin{equation}
  \label{bon}
  df_{\pm} = \overline{(f^{*} \pm \epsilon)} df (f^{*} \pm \epsilon), 
\end{equation}
where $\epsilon > 0$ is a real constant, $f$ is a global isothermic 
immersion defined on the universal cover of $M$ and $f^{*}$  is 
dual to $f$ \cite{pipeka}. The dual $f^{*}$ descends to a map on $M$. 

The difference between the second fundamental forms of two Bonnet
mates $f_{\pm}$ defines a holomorphic quadratic differential, called
the shape distortion differential $D = \left(\II_{+} -
  \II_{-}\right)^{2,0}$.  The zeroes of $D$ are precisely the
umbilic points of $f_{\pm}$ (see \cite{kamberov-bonnet-96}).
\begin{observation} \label{shapedist} 
Suppose that $f_{\pm}$ are Bonnet mates generated by the isothermic 
surface $f$, that is $df_{\pm} = \overline{\lambda_{\pm}} df \lambda_{\pm}$, 
where $\lambda_{\pm} = f^{*} \pm \epsilon$, and $\epsilon$ is a
positive constant. Then 
\begin{equation} \label{shaped}
df\backslash D = 4\epsilon*df^{*}
\end{equation}
\end{observation} 
{\bf Proof} \, Remark~\ref{hopf} implies  $df\backslash D = df\backslash -\nu_{f}^{-1}\circ 
\left(
\nu_{f_{+}}\circ df_{+} \omega_{+} - \nu_{f_{-}}\circ df_{-} \omega_{-}
\right)$. 
Therefore, from the following lemma:
\begin{lemma} \label{lemma:gauge} Suppose that $\tilde{f}$ is obtained from $f$ by a spin transform, thus
  $d\tilde{f} = \overline{\lambda}df\lambda$. Let $\tilde{\tau}$ be tangential 
and anti-conformal with respect to $\tilde{f}$, then 
 $ df\backslash \nu^{-1}_{f}\circ \nu_{\tilde{f}} \circ d\tilde{f} \tilde{\tau} = 
\lambda \tilde{\tau}\overline{\lambda}.$
\end{lemma}
it follows that
\begin{equation} \label{eq:differ}
df\backslash D = -\left(\lambda_{+} \omega_{+} \overline{\lambda_{+}} -  \lambda_{-} \omega_{-} \overline{\lambda_{-}}\right).
\end{equation}
Let $N_{\pm}$ denote the normal to $f_{\pm}$, and so $N = \lambda_{\pm} N_{\pm} \lambda_{\pm}^{-1}$.
 From \mref{eqn:Weingarten} it follows
\begin{equation}
  \label{eq:2}
  \lambda_{\pm} \omega_{\pm} \overline{\lambda_{\pm}}  = \lambda_{\pm} dN_{\pm} \overline{\lambda_{\pm}} + 
H_{\pm} \lambda_{\pm} df_{\pm} \overline{\lambda_{\pm}}. 
\end{equation}
But $df_{\pm} = \overline{\lambda_{\pm}} df \lambda_{\pm}$, and so 
$H_{\pm} \lambda_{\pm} df_{\pm} \overline{\lambda_{\pm}} = H_{\pm}|\lambda_{\pm}|^{4} df$. The immersions 
$f_{\pm}$ are Bonnet mates and hence $H_{+}=H_{-}$ and $|\lambda_{+}| = |\lambda_{-}|$. Thus \mref{eq:2} implies 
\begin{equation}
  \label{eq:3}
  \lambda_{+} \omega_{+} \overline{\lambda_{+}} -  \lambda_{-} \omega_{-} \overline{\lambda_{-}} = 
  \lambda_{+} dN_{+} \overline{\lambda_{+}} -  \lambda_{-} dN_{-} \overline{\lambda_{-}}.
\end{equation}
Using, $d\lambda^{-1} = -\lambda^{-1} d\lambda \lambda^{-1}$ we get 
\begin{equation} \label{eq:4}
\lambda_{\pm} dN_{\pm} \overline{\lambda_{\pm}} = 
\lambda_{\pm} d(\lambda_{\pm}^{-1} N \lambda_{\pm} ) \overline{\lambda_{\pm} } \nonumber \\
  = \mbox{} - (d\lambda_{\pm}) \lambda_{\pm}^{-1} N |\lambda_{\pm}|^{2} + dN |\lambda_{\pm}|^{2} + 
N(d\lambda_{\pm})\overline{\lambda_{\pm} }.
\end{equation}
Recall that $\lambda_{\pm} = f^{*} \pm \epsilon$ and so $d\lambda{\pm} = df^{*}$ and hence 
\[
\lambda_{+} \omega_{+} \overline{\lambda_{+}} -  \lambda_{-} \omega_{-} \overline{\lambda_{-}} = 
-d\lambda_{+} 2\epsilon N + Ndf^{*} 2\epsilon = -4\epsilon *df^{*},
\]
and so $df\backslash D = 4\epsilon *df^{*}$. Finally we prove the lemma. \\
{\bf Proof of Lemma}~\ref{lemma:gauge} 
Let $a$ and $b$ the real valued one-forms on $M$ such that 
$ d\tilde{f} \, \tilde{\tau} = a + b\tilde{N}$, then 
\[
\nu^{-1}_{f}\circ \nu_{\tilde{f}} \circ d\tilde{f} \, \tilde{\tau} = \nu^{-1}_{f} (a + ib) = a + bN.
\]
Substituting $d\tilde{f} \, \tilde{\tau} = \lambda^{-1} (a + bN)\lambda$ in  
$
df \, \lambda\tilde{\tau}\overline{\lambda} = (\overline{\lambda}^{-1}) d\tilde{f} \lambda^{-1}\lambda\tilde{\tau}\overline{\lambda}= 
\frac{1}{|\lambda|^2} (a + bN) |\lambda|^{2}
$ 
yields $df \, \lambda\tilde{\tau}\overline{\lambda} = a + bN = \nu^{-1}_{f}\circ \nu_{\tilde{f}} \circ d\tilde{f} \, \tilde{\tau}$.
\QED 
\begin{remark} \label{branchpoints} \, Observation~\ref{shapedist} 
implies that the branch points of $f^*$ are precisely the zeroes of the
shape distortion operator $D$ and therefore they are precisely the 
umbilic points of $f_{\pm}$ (\cite{kamberov-bonnet-96}). The  
branch points of $f^*$ are umbilic points of $f$ but $f$ could have other umbilic 
points too.  Moreover, the order of
vanishing of $\omega$ at the branch points is at least equal to the
order of vanishing of $df^*$, in particular it follows that $\omega=adf^*$, 
where the real valued coefficient $a$ is bounded, possibly continuous, 
at the branch points of $f^*$, and smooth away from them.
\end{remark}
Now we can give the construction which yields all Bonnet 
surfaces which are not included in the Bonnet-Cartan-Chern classification. 
By necessity these surfaces must have isolated umbilics and their mean 
curvature is not constant in a neighborhood of the umbilic. 
\begin{theorem} \label{umbconstr}  \ \\

{\bf (i.)} Let $f_{\pm}$ be two Bonnet mates defined by \mref{bon}. A 
 point $p_{0}$ is umbilic for $f_{\pm}$ if and only
if $p_0$ is a branch point for $f^*$. \\

{\bf (ii.)} If $f$ is a non-minimal isothermic immersion with a dual $f^*$ then 
there exists at most one $\epsilon > 0$ such that the Bonnet mates $f_\pm$ 
  defined by \mref{bon} have constant mean curvature. 
In particular, let $f$ be an immersion with positive constant mean curvature
and let $f^{*}$ be a dual map for $f$, then for every $\epsilon$ the Bonnet
mates $f_{\pm}$ defined by \mref{bon} have non-constant mean curvature. 
\end{theorem}
{\bf Proof} \,  Part {\bf (i)} follows from Remark~\ref{branchpoints}. 
To prove part {\bf (ii)} fix a non-minimal isothermic immersion $f$, and a dual
$f^{*}$.  If there exists a positive $\epsilon$ such that the
corresponding Bonnet mates $f_{\pm}$ have constant mean curvature $c$, then 
$c=H/(|f^{*}|^2 + \epsilon^2)$, and so $c$ is the unique solution of
$dH=cd|f^{*}|^{2}$ (if such exists). But then $\epsilon$ is uniquely
determined also.  \QED

Thus all new Bonnet surfaces are obtained by applying the \cite{pipeka}
construction to isothermic immersions whose duals have branch points. 
To produce explicit examples simply take  $f$ to be the Mr. Bubble constant mean
curvature immersion. 
\begin{coro} There exist non-constant mean curvature Bonnet mates with 
umbilic points. 
\end{coro}
\begin{remark} Since isothermic surfaces could be only $C^{1}$,  
the construction in \cite{pipeka} indicates that we can make sense of ``Bonnet surfaces'' 
with very low smoothness, say $C^{1}$. The classical results required at least $C^{5}$. 
\cite{Chern-Deform-Pr-Curv}. (See \cite{sqd} for further discussion.)
\end{remark}
\bigskip
\noindent
{\bf Acknowledgments}

\noindent
{\footnotesize
We thank Sasha Bobenko, Gary Jensen,  Udo Hetrich-Jeromin, Gerda Kamberov, 
Peter Norman, Franz Pedit, Ulrich Pinkall for their interest and for many discussions, 
and to Robert Bryant for introducing us to the Bonnet 
problem and for his interest. 
}
\bigskip
\noindent

\bibliographystyle{plain}

\end{document}